\documentstyle[amssymb,righttag]{amsart}
\numberwithin{equation}{section}
\numberwithin{figure}{section}
\newtheorem{theorem}{Theorem}[section]
\newtheorem{lemma}[theorem]{Lemma}

\theoremstyle{definition}

\theoremstyle{remark}

\numberwithin{equation}{section}

\begin{document}

\title[The Markov-Bernstein inequalities in $ L^{2} \left((-1, 1), \,w^{(\alpha,\;\beta)}\right) $. ] {Asymptotics of sharp constants \\of
Markov-Bernstein inequalities\\ in  integral norm  with Jacobi weight}

\thanks{The  first and fourth authors were partly supported by the program  N1 of DMS RAS
and grant RFBR-11-01-00245. The third author was partly supported by  the Scientific Schools program - 4664.2012.1.  The paper has been finished while the first author visited INSA, Rouen, France.}

\author[A. Aptekarev]{A. I. Aptekarev }

\author[A. Draux]{A. Draux}

\author[V. Kalyagin]{V. A. Kalyagin}

\author[D. Tulyakov]{D. N.  Tulyakov}




\begin{abstract} The classical A. Markov inequality establishes a relation between  the maximum modulus
or  the $L^{\infty}\left([-1,1]\right)$ norm of a polynomial $Q_{n}$ and  of its derivative:
 $\|Q'_{n}\|\leqslant M_{n} n^{2}\|Q_{n}\|$, where the constant $M_{n}=1$ is sharp.
 The limiting behavior of the sharp constants $M_{n}$
 for this inequality, considered in the space $L^{2}\left([-1,1], w^{(\alpha,\beta)}\right)$ with respect to
 the classical Jacobi weight
 $w^{(\alpha,\beta)}(x):=(1-x)^{\alpha}(x+1)^{\beta}$,
 is studied. We prove that, under the condition $|\alpha - \beta| < 4 $, the limit is
 $\lim_{n \to \infty} M_{n} = 1/(2 j_{\nu})$ where $j_{\nu}$ is
    the smallest
 zero of the Bessel function $J_{\nu}(x)$ and
 $2 \nu= \mbox{min}(\alpha, \beta) - 1$.
\end{abstract}

\maketitle

\section{Introduction}\label{sec:1}
A quantity
\begin{equation}\label{1}
M_{n}:=\sup\limits_{\deg\,Q_{n}\leqslant n}\frac{\|Q'_{n}\|_{\mathfrak{X}_{1}}}{\|Q_{n}\|_{\mathfrak{X}_{2}}}\;,
\quad Q_{n}\; \mbox{- polynomial},
\end{equation}
is called the sharp constant for the Markov-Bernstein inequality in
functional spaces $\mathfrak{X}_{1},\mathfrak{X}_{2}$ with norms
$\|.\|_{\mathfrak{X}_{1}}$, $\|.\|_{\mathfrak{X}_{2}}$.

The remarkable classical inequality of A.A. Markov
for $\mathfrak{X}_{1}=\mathfrak{X}_{2}=L^{\infty}\left([-1,1]\right)$
$$\|Q'_{n}\|\leqslant \, n^{2}\|Q_{n}\|$$ is sharp \cite{Markov}.
We recall that the corresponding inequality for the trigonometric polynomials  has firstly been obtained by S.N. Bernshtein in \cite{Bernshtein}.
His result was not sharp, and the sharp version is due to E. Landau (see \cite{Akhiezer}). For the weighted $L^2$ spaces $\mathfrak{X}_{1}=\mathfrak{X}_{2}:=L^{2}\left([a, b], w\right)$, for some classical weights, the sharp constants \eqref{1} are known (see \cite{MMR}):
\begin{itemize}
\item [\textbf{1.}\,] $w = \exp(-x^2),\,\, x\in (-\infty, \infty), \,\, M_n =\sqrt{2n}\qquad \qquad $ (E. Schmidt, 1944);
    \item [\textbf{2.}\,] $w = \exp(-x),\,\,\,\, x\in (0, \infty), \,\,\,\, M_n \,=\, \displaystyle\frac{1}{2\sin(\frac{\pi}{4n+2})} \,\,\, $ (P. Turan, 1960).
\end{itemize}
However, for other classical weights, explicit expressions for the sharp constants are not known. In \cite{ADK} results on the asymptotics
$
M_{n}\rightarrow \,\,?\,,\; n\rightarrow\infty\;
$
were discussed. In particularly, for the Gegenbauer weight $w^{(\alpha)}(x):=(1-x^2)^{\alpha}\;,\; x\in[-1,1],  \alpha >-1,
$ the following result was stated there
\begin{equation}\label{Gegen}
M_{n}=\frac{n^2}{2 j_{\,\nu (\alpha)}}\;\left(1\,+\, o\; (1)\right)\,,\qquad
\nu (\alpha)\,:= \,\frac{\alpha - 1}{2}\,,
\end{equation}
where  $j_{\,\nu}$ is the smallest
 zero of the Bessel function $J_{\nu}(x)$ (we shall keep the notations $\nu(\alpha)$   and $ j_{\,\nu}$ in what follows).

In the present paper we study the asymptotics of the sharp constant \eqref{1}  for the classical Jacobi weight (the space $L^{2}([-1, 1], w^{(\alpha,\,\beta)}$ is defined in the section \ref{subsec:2}):
\begin{equation}\label{JacW}
\mathfrak{X}_{1}=\mathfrak{X}_{2}:=L^{2}\left([-1, 1], w^{(\alpha,\,\beta)}\right),\; w^{(\alpha,\,\beta)}(x):=(1-x)^{\alpha}(x+1)^{\beta}, \;\alpha,\beta>-1.
\end{equation}
The main result of our paper is
\begin{theorem}\label{T1} Let the parameters of the Jacobi weight \eqref{JacW} satisfy the restriction
\begin{equation}\label{JacWrestr}
w^{(\alpha,\,\beta)}(x)\,\,:\qquad |\alpha - \beta| < 4.
\end{equation}
Then, for the sharp constant \eqref{1} in the space \eqref{JacW}, we have the asymptotics:
\begin{equation}\label{AsSCJac}
M_{n}=\frac{n^2}{2 j_{\,\nu^*}}\;\left(1\,+\, o\; (1)\right)\,,\qquad
\nu^*= \min \{\nu (\alpha), \,\nu(\beta)\}\,.
\end{equation}
\end{theorem} \noindent
We see that, for $\alpha = \beta$, asymptotics \eqref{AsSCJac} match \eqref{Gegen}. When $\alpha \ne \beta$, then asymptotics \eqref{AsSCJac} look as a reasonable generalization of \eqref{Gegen}. The most surprisingly for us is the appearance of the restriction  \eqref{JacWrestr}. At the moment we can not prove or disapprove its necessity, however, we have  to admit that this restriction is unavoidable in our proof strategy of  Theorem~\ref{T1}.

The rest of the paper contains the proof of  Theorem~\ref{T1}. Our approach consists on the following steps:
\begin{itemize}
\item [\textbf{1)}\,] We start with an explicit representation of  $M_n$ as the eigenvalue of a linear operator in $\mathbb{R}^{n}$ defined by a five diagonal matrix.
\item [\textbf{2)}\,] Then we state a Finite Difference (FD) Boundary Value Problem (BVP) which is equivalent to the  eigenvalue problem.
\item [\textbf{3)}\,] The next step is to determine a limiting (for FD problem) Differential Equation (DE) and  its general solution.
\item [\textbf{4)}\,] Then we vanish the spectral parameter in  FD problem and find  linearly independent Particular Solutions (PS) satisfying Boundary Conditions (BC) at the initial values of the discrete variable (the left end BC). For the small (with respect to $n$) indexes in FD and spectral parameter in the fixed range the asymptotics of the solutions of FD does not depend on the spectral parameter. Therefore the initial conditions can be rewritten as asymptotics condition for the indexes $1 << k << n$. Then this condition is exported to the boundary condition of DE.
\item [\textbf{5)}\,] Matching these  FD problems, we get PS of the limiting DE.
\item [\textbf{6)}\,] Finally, taking these PS of DE as an approximation of the PS of FD, we satisfy the right end BC of FD BVP. It gives an approximation of the desired eigenvalue like in  \eqref{Gegen} or \eqref{AsSCJac}.
\end{itemize}

These steps are performed in the next Section~\ref{sec:2}. Some of these steps have already been studied before for various functional spaces  in \eqref{1}, see in \cite{MMR} Chapter 6 Section 6.1.6, \cite{DE}, \cite{ADK}, \cite{ATD}, \cite{DK}.

However, to conclude a rigorous proof of the Theorem~\ref{T1}, it remains to justify the final step, i.e. to prove that  PS of DE\ which match the satisfying to the left end BC of
FD problem, indeed are close  to PS of FD problem. In Section~\ref{sec:3} we state and prove the corresponding result, see Theorem~\ref{Th2}. This theorem establishes a new result on the local asymptotics of the powerlike growing solution of the high order recurrence relations. Previous results in this direction are in \cite{Apt}, \cite{Tul}, \cite{Tul1}. In what follows we use the notation $(a)!$ for the value of $\Gamma(a+1)$.

\section{Finite difference BVP for $M_n$ and its differential approximation}\label{sec:2}
\subsection{A spectral representation for $M_{n}$ in $\mathbb{R}^{n}$.}\label{subsec:2}
We note from \eqref{1}, that $M_n$ is the norm of the operator differentiation in a finite dimensional space   $\mathcal{P}_{n}$ of polynomials of degree at most equal to $n$.
Let $Q_{n}$ be an arbitrary polynomial of $\mathcal{P}_{n}$. We take the expansion of this polynomial $Q_{n}$ and of its derivative $Q'_{n}$ in the basis of monic Jacobi polynomials $P_{k}^{(\alpha,\beta)}$ ($P_{k}^{(\alpha,\beta)}(x)=x^{k}+\cdots$).
By using the $L^{2}_{w^{\alpha,\beta}}$ inner product
\begin{equation}\label{2}
(g,f):=\int\limits_{-1}^{1}g(x)\,\overline{f(x)}\,w^{(\alpha,\beta)}(x)\,dx\;,
\end{equation}
the square norm of $P_{n}^{(\alpha,\beta)}$ is
\begin{equation}\label{4w}
\|P_{n}^{(\alpha,\beta)}\|^{2}=(P_{n}^{(\alpha,\beta)},P_{n}^{(\alpha,\beta)} )=2^{2n+\alpha+\beta-1}\displaystyle\frac{n!(n+\alpha)!
(n+\beta)!(n+\alpha+\beta)!}{(2n+\alpha+\beta)!(2n+\alpha+\beta+1)!}.
\end{equation}
Then, we have (in general $Q_n=\sum_{k=0}^n c_kP_k^{(\alpha,\beta)}$ but to solve (\ref{1}) it is sufficient to consider the case where $c_0=0$):
\begin{equation}\label{5w}
Q'_{n}(x):=\sum\limits_{k=0}^{n-1}v_{k}P_{k}^{(\alpha,\beta)}\;,\quad
Q_{n}(x):=\sum\limits_{k=0}^{n-1}u_{k}P_{k+1}^{(\alpha,\beta)}\;.
\end{equation}
Differentiating $Q_{n}$ here, and using the property of Jacobi
polynomials
$$
\frac{d}{dx}P_{k}^{(\alpha,\beta)}(x)=kP_{k-1}^{(\alpha+1,\beta+1)}(x)\;,
$$
we arrive to
$$
Q'_{n}=\sum\limits_{k=0}^{n-1}u_{k}P_{k+1}^{'(\alpha,\beta)}=
\sum\limits_{k=0}^{n-1}(k+1)u_{k}P_{k}^{(\alpha+1,\beta+1)}
=\sum\limits_{k=0}^{n-1}v_{k}P_{k}^{(\alpha,\beta)}\;.
$$
Then, applying  the well known $\alpha$ increasing (and $\beta$ increasing) relation (see \cite{ABST}, Chapter 22):
\begin{equation}\label{6w}
  P_{n}^{(\alpha,\beta)}=P_{n}^{(\alpha+1,\beta)}-\displaystyle\frac{2n(n+\beta)P_{n-1}^{(\alpha+1,\beta)}}
  {(2n+\alpha+\beta)(2n+\alpha+\beta+1)}\,,
\end{equation}
we obtain for the vectors from \eqref{5w}
$$
\vec{v}:=(v_{0},v_{1},\ldots,v_{n-1})^{Tr}\;,\quad\vec{u}:=(u_{0},u_{1},\ldots,u_{n-1})^{Tr}\;,
$$
the following  relations
\begin{equation}\label{7}
\mathbf{N}\vec{u}=\mathbf{C}_{2}\mathbf{C}_{1}\vec{v}\;,\qquad \qquad \mathbf{N}:=\mbox{diag}(1,2,\ldots,n)\;,
\end{equation}
and for the $n \times n$ matrices $\mathbf{C}_{2},\mathbf{C}_{1}$ we have from (\ref{6w})
\begin{equation*}
\begin{array}{cc}
  \begin{array}{l}
     \mathbf{C}_{1}:=\mathbf{I}-\mbox{diag}\left(\frac{2k(k+\beta)}{(2k+\alpha+\beta)
     (2k+\alpha+\beta+1)}\right)_{k=1}^{n} \mathbf{T }\\
     \\
     \mathbf{C}_{2}:=\mathbf{I}+\mbox{diag}\left(\frac{2k(k+\alpha+1)}{(2k+\alpha+\beta+1)
     (2k+\alpha+\beta+2)}\right)_{k=1}^{n} \mathbf{T} \\
     \end{array} ,
   & \mathbf{T}:=\left(\begin{array}{ccccc}
                0 & 1 &  &  & 0 \\
                 & 0 & 1 &  &  \\
                 &  & 0 & \ddots &  \\
                0 &  & & \ddots & 1 \\
                 0 & 0 &  &  & 0
              \end{array}
   \right)
\end{array}.
\end{equation*}
Now, we write norms for (\ref{5w}) by using the inner product  (\ref{2}):
$$
\begin{array}{cc}
  \|Q'_{n}\|^{2}=<\vec{v},\mathbf{D}\vec{v}>\;, & \mathbf{D}:=\mbox{diag}\;\left(\|P_{k}^{(\alpha,\beta)}\|^{2}\right)_{k=0}^{n-1}\;, \\
\\
  \|Q_{n}\|^{2}=<\vec{u},\mathbf{D}^{+}\vec{u}>\;, & \mathbf{D}^{+}:=\mbox{diag}\;\left(\|P_{k}^{(\alpha,\beta)}\|^{2}\right)_{k=1}^{n}\;.
\end{array}
$$
Where $<\ ,\ >$ is usual vector inner product. Thus,  for the sharp constant in (\ref{1}) - \eqref{JacW} we have by using (\ref{7})
\begin{equation}\label{9w}
M_{n}^{2}=\sup\limits_{\vec{v}}\frac{<\vec{v},\mathbf{D}\vec{v}>}{<\mathbf{N}^{-1}
\mathbf{C}_{2}\mathbf{C}_{1}\vec{v},
\mathbf{D}^{+}\mathbf{N}^{-1}\mathbf{C}_{2}\mathbf{C}_{1}\vec{v}>}\;=\;
\sup\limits_{\vec{v}}\frac{<\vec{v},\mathbf{D}\vec{v}>}{<\vec{v},\mathbf{A}\vec{v}>}\;,
\end{equation}
where we denote
\begin{equation}\label{10}
\mathbf{A}:=\mathbf{C}_{1}^{Tr}\mathbf{C}_{2}^{Tr}\mathbf{N}^{-1}\mathbf{D}^{+}\mathbf{N}^{-1}\mathbf{C}_{2}\mathbf{C}_{1}
,\,\,\,
\mathbf{D}:=\mbox{diag}(d_{k})_{k=0}^{n-1},\, \mathbf{D}^{+}:=\mbox{diag}(d_{k})_{k=1}^{n}.
\end{equation}
For the purpose of (\ref{9w}) we can omit the factor $2^{\alpha+\beta-1}$ in
(\ref{4w}), so we put
\begin{equation}\label{11}
d_{k}:=\frac{2^{2k}k!(k+\alpha)!(k+\beta)!(k+\alpha+\beta)!}{(2k+\alpha+\beta)!(2k+\alpha+\beta+1)!}
\,.
\end{equation}
Finally, from (\ref{9w}) we get by the arguments of pencil of quadratic forms (see \cite{GANT}, Chapter 10.7) the spectral radius representation for
the exact constant:
\begin{equation}\label{12}
M_{n}^{2}=\lambda_{\min}^{-1}(\mathbf{A},\mathbf{D})\;,
\end{equation}
where $\lambda_{\min}(\mathbf{A},\mathbf{D})$ is a root (with the minimal modulus)  of
the equation
\begin{equation}\label{13}
\det(\mathbf{A}-\lambda \mathbf{D})=0\;,
\end{equation}
and correspondingly the eigenvector $\vec{v}$
\begin{equation}\label{14}
(\mathbf{A}-\lambda_{\min}\mathbf{D})\vec{v}=0
\end{equation}
defines the extremal polynomial (\ref{5w}) in (\ref{1}), (\ref{JacW} ).

\subsection{Finite difference equation for the coordinates of $\vec{v}$.}\label{subsec:3}
To simplify expressions (i.e. to cancel factorials) in what follows, we introduce a new variable
for the coordinates of the vector $\vec{v}$ (see (\ref{5w}):
\begin{equation}\label{15}
v_{k}=:x_{k}\,\displaystyle\frac{(2k+\alpha+\beta+1)!}{2^{k}(k+\alpha)!(k+\beta)!}\;,\qquad k=0,\ldots,n-1\;.
\end{equation}
Next, taking the $k$-th coordinate of the equation (\ref{14})
\begin{equation}\label{5-term}
\left[\,(\mathbf{A}-\lambda \mathbf{D})\;\vec{v}\,\right]_{k}=0\;, \qquad k=0,\ldots,n-1\,
\end{equation}
we get a 5-term recurrence relation which connects the coordinates
$\{x_s\}^{k+2}_{s=k-2}$, $k=0,1,2,\ldots,n-1$:
\begin{equation*}
  \begin{array}{l}
x_{k+2}\frac{(k+2)!}{(k-2)!}\frac{(2k+\alpha+\beta)!}{(2k+\alpha+\beta-3)!}
\frac{(k+\alpha+\beta+1)!}{(k+\alpha+\beta-1)!}\,=\,
x_{k+1}\frac{(k+1)!}{(k-2)!}\frac{(2k+\alpha+\beta-1)!}{(2k+\alpha+\beta-3)!}\,\Xi_{1}
\,+\\\\
x_{k}\left(\Xi_{2}-\lambda\,\frac{(k+1)!}{(k-2)!}\frac{(2k+\alpha+\beta+4)!}{(2k+\alpha+\beta-3)!}
\frac{(k+\alpha+\beta)}{4}\right)+
x_{k-2}\frac{(k+1)!}{(k-1)!}\frac{(2k+\alpha+\beta+4)!}{(2k+\alpha+\beta+1)!}
\frac{(k+\alpha)!}{(k+\alpha-2)!}\frac{(k+\beta)!}{(k+\beta-2)!}\\\\
+\quad x_{k-1}\frac{(2k+\alpha+\beta+4)!}{(2k+\alpha+\beta+2)!}\,
(k^2-1)(k+\alpha)(k+\beta)(2k+\alpha+\beta-1)(\alpha+\beta-2)(\alpha-\beta)\,,
\end{array}
\end{equation*}
where
\begin{equation*}
  \begin{array}{l}
\Xi_{1}= (k+\alpha+\beta)(2k+\alpha+\beta+3)(\alpha+\beta-2)(\alpha-\beta),\,\,
\Xi_{2}= \displaystyle\frac{k^4}{2}+ k^3(1+\alpha+\beta)
 \\\\ +\,\, k^2\,\,\frac{2\alpha+2\beta+2\alpha^2+3\alpha \beta+2\beta^2+1}{2}\,\, + \,\,k\,\, \frac{(1+\alpha+\beta)(\alpha^2+\alpha \beta+\beta^2)}{2}
 \,\,+ \, \,O_{\alpha, \;\beta}(1).
\end{array}
\end{equation*}
These finite difference equation can be considered as spectral equation for the problem (\ref{14}). We obtain a non trivial solution of the (\ref{14}) if $x_{-1}=x_{-2}=0$ and $x_{n}=x_{n+1}=0$. These boundary conditions will be widely used in the paper.
The 5-terms recurrence equation can be rewritten  in a matrix form for the bundle $\overrightarrow{X}_{k}$:
\begin{equation}\label{16}
\overrightarrow{X}_{k+2}=\mathbf{M}_{2}(k,\lambda)\overrightarrow{X}_{k}\;,\quad \overrightarrow{X}_{k}:=
\left(
x_{k-2},\,
x_{k-1},\,
x_{k},\,
x_{k+1}\,
\right)^{Tr},\quad k\in\mathbb{Z}_{+}\;.
\end{equation}
The matrix $\mathbf{M}_{2}$, can be divided in two terms (one is linearly dependent on $\lambda$, the other one is independent on $\lambda$) $\mathbf{M}_{2}(k,\lambda)=\lambda \mathbf{M}_{2}^{(1)}+\mathbf{M}_{2}^{(0)}$.

The leading coefficients of the expansion of
the matrices $\mathbf{M}_{2}^{(1)}$ and $\mathbf{M}_{2}^{(0)}$ are
$$
\mathbf{M}_{2}^{(1)}:=\left(
               \begin{array}{cccc}
                 0 & 0 & 0 & 0 \\
                 0 & 0 & 0 & 0 \\
                 0 & 0 & (\tilde{b}-8k-2a-3) & 0 \\
                 0 & 0 & 4(2-a)(\alpha-\beta) & (\tilde{b}-16k-6a-15) \\
               \end{array}
             \right)+O\left(\frac{1}{k}\right)\;,
$$
where $a:=\alpha+\beta, \quad \tilde{b}:=-4k^{2}-4ak-2a^{2},\,\,$ and
$$
\mathbf{M}_{2}^{(0)}:=\left(
               \begin{array}{cccc}
                 0 & 0 & 1 & 0 \\
                 0 & 0 & 0 & 1 \\
                 -1-\frac{2}{k} & 0 & 2+\frac{2}{k} & 0 \\
                 0 & -1-\frac{2}{k} & 0 & 2+\frac{2}{k} \\
               \end{array}
             \right)+O\left(\frac{1}{k^{2}}\right)\;.
$$
The  relations (\ref{16}) can be rewritten as a finite
difference equation  involving the vectors $\overrightarrow{X}_{k}$:
$$
\displaystyle\frac{\overrightarrow{X}_{k+2}-\overrightarrow{X}_{k}}{2/n}=\frac{n}{2}[\mathbf{M}_{2}
-\mathbf{I}]\overrightarrow{X}_{k}\;,\qquad k=0,2,\ldots\;\;.
$$
In order to work with a better structured matrix,we pass from $\overrightarrow{X}_{k}$ to
$\overrightarrow{Y}_{k}$:
\begin{equation}\label{17}
\overrightarrow{Y}_{k}:=\mathbf{U}_{k}\overrightarrow{X}_{k}\;,\quad \mathbf{U}_{k}:=\left(
 \begin{array}{cccc}
 1 & 1 & 0 & 0 \\
 1 & -1 & 0 & 0 \\
 -k & -k & k & k \\
 -k & k & k & -k \\
 \end{array}
 \right)\ \ k>1,
\quad \mathbf{U}_{0}:=\left(
 \begin{array}{cccc}
 1 & 1 & 0 & 0 \\
 1 & -1 & 0 & 0 \\
 -1 & -1 & 1 & 1 \\
 -1 & 1 & 1 & -1 \\
 \end{array}
 \right)
\end{equation}
satisfying
\begin{equation}\label{17a}
\overrightarrow{Y}_{k+2}=\widehat{\mathbf{M}}_{2}^{(\alpha,\beta)}(k,\lambda)\overrightarrow{Y}_{k}\;,\qquad
\widehat{\mathbf{M}}_{2}^{(\alpha,\beta)}:=\mathbf{U}_{k+2}\mathbf{M}_{2}\mathbf{U}_{k}^{-1}.
\end{equation}
Then, we arrive to a finite-difference system $k=2,4,\ldots$,
\begin{equation}\label{18}
\displaystyle\frac{\overrightarrow{Y}_{k+2}-\overrightarrow{Y}_{k}}{2/n}=\frac{n}{k}
\mathbf{M}_{3}(k,\lambda)\overrightarrow{Y}_{k}\;,\quad
\overrightarrow{Y}_{2}=\overrightarrow{Y}_{0}+2\mathbf{M}_3(0,\lambda)\overrightarrow{Y}_{0},
\quad \overrightarrow{Y}_{0}=(0,0,C_1,C_2)^{Tr}
\end{equation}
where $\mathbf{M}_3(0, \lambda)=(1/2)[U_2\mathbf{M}_2(0,\lambda)U_0^{-1}-\mathbf{I}]$ and the matrix
$
\mathbf{M}_{3}(k)=\displaystyle\frac{k}{2}\,[\mathbf{U}_{k+2}\mathbf{M}_{2}\mathbf{U}_{k}^{-1}-\mathbf{I} ]=\lambda \mathbf{M}_{3}^{(1)}(k)+\mathbf{M}_{3}^{(0)}(k)
$
 ($k>0$) has expansions  $\mathbf{M}_{3}^{(0)}:=$
$$
\left(
\begin{array}{cccc}
0 & 0 & 1/2 & 0 \\
0 & 0 & 0 & 1/2 \\
2\alpha(\alpha-2)-\frac{\alpha\,\diamondsuit(\alpha,\beta)}{k}
& \frac{2\beta(\beta-2)}{k} & 2+\frac{\Box(\alpha,\beta)}{2k} & \frac{\triangle(\alpha,\beta)}{2k} \\\\
\frac{2\alpha(\alpha-2)}{k} & 2\beta(\beta-2)-\frac{\beta\,\diamondsuit(\beta,\alpha)}{k}
& \frac{\triangle(\beta,\alpha)}{2k} & 2+\frac{\Box(\beta,\alpha)}{2k} \\
\end{array}
\right)\,+\,
O\left(\frac{1}{k^{2}}\right)\;,
$$
and
$$
\mathbf{M}_{3}^{(1)}:=\left(
\begin{array}{cccc}
0 & 0 & 0 & 0 \\
0 & 0 & 0 & 0 \\
-2k^{4}-(2a+10)k^{3} & 2k^{3} & -2k^{3} &  2k^{2} \\
 2k^{3} & -2k^{4}-(2a+10)k^{3} & 2k^{2} & -2k^{3} \\
\end{array}
\right)+O(k)\;.
$$
Here we denoted
$$
  \diamondsuit(\alpha,\beta)=4-6\alpha+2\alpha^{2}+2\alpha\beta-\beta\,, \qquad
  \Box(\alpha,\beta)=4\alpha^{2}-9\alpha+2\alpha\beta-\beta+4 \,,
$$$$
  \triangle(\alpha,\beta)=\alpha^{2}-\beta^{2}-2\alpha+2\beta+1  \;
, \qquad \text{and} \qquad a\,:=\,\alpha+\beta\,.\qquad \qquad \qquad \quad
$$
\subsection{General solution of the limiting system of ODEs.}\label{4}
Now, we take a formal limit (under an appropriate
scalling) of the Finite Difference (FD) problem (\ref{18}) to
arrive to a limiting system of  ordinary differential equations
(ODEs). Indeed, if we denote
$\vec{y}=(y_{1},y_{2},y_{3},y_{4})^{Tr}$:
\begin{equation}\label{20}
\vec{y}(t,l):=\lim\limits_{n\to\infty,\,\frac{k}{n}\to t}\overrightarrow{Y}_{k}(\lambda)\Bigr|_
{\lambda=l / n^{4}}\;,
\end{equation}
(we shall investigate the existence of this limit later),
then we arrive from (\ref{18}) to the system of ODEs:
\begin{equation}\label{21}
\frac{d}{dt}\,\vec{y}(t,l)=\frac{1}{t}\widetilde{\mathbf{M}}_{3}(t,l)\,\vec{y}(t,l),\quad
\widetilde{\mathbf{M}}_{3}(t,l)=\lim\limits_{\frac{k}{n}\to t}\mathbf{M}_{3}(k,\frac{l}{n^{4}})\;,
\end{equation}
where
$$
\widetilde{\mathbf{M}}_{3}(t,l)=\left(
\begin{array}{cccc}
0 & 0 & 1/2 & 0 \\
0 & 0 & 0 & 1/2 \\
-2lt^{4}+2\alpha(\alpha-2) & 0 & 2 & 0 \\
0 & -2lt^{4}+2\beta(\beta-2) & 0 & 2 \\
\end{array}
\right).
$$
Due to the special structure of the main term of asymptotics of the matrix $M_{3}$ in (\ref{18}),
this system is split in two second order independent scalar ODEs:
$$
\left\{\begin{array}{l}
t\,\displaystyle\frac{d}{dt}\,y_{j}(t)=\displaystyle\frac{z_{j}(t)}{2}\;, \\
\\
t\,\displaystyle\frac{d}{dt}\,z_{j}(t)=2[(b_{j}(b_{j}-2)-lt^{4})y_{j}(t)+z_{j}(t)]\;,
\end{array}
\right. \quad j=1,2.
$$
and
\begin{equation}\label{22}
\displaystyle\frac{d^{2}}{dt^{2}}\,y_{j}(t,l)=\frac{1}{t}\displaystyle\frac{d}{dt}\,y_{j}(t,l)-
\left(t^{2}l-\frac{b_{j}(b_{j}-2)}{t^2}\right)\,y_{j}\;,\quad j=1,2.
\end{equation}
Here and in what follows, we use the notation:
\begin{equation}\label{notbj}
b_{j}:=\left\{\begin{array}{c}
                \alpha,\;j=1 \\
                \beta,\;j=2
              \end{array}
\right.\;.
\end{equation}
The ODE (\ref{22}) is a modified Bessel equation in Bowman form (see \cite{BOW}, Ch. 104); its general solution is
\begin{equation}\label{23}
y_{j}(t,l)=\widetilde{C}_{1,j}\,t\,J_{\nu(b_{j})}\left(\sqrt{l}\,\,\frac{t^{2}}{2}\right)+
\widetilde{C}_{2,j}\,t\,Y_{\nu(b_{j})}\left(\sqrt{l}\,\,\frac{t^{2}}{2}\right)\;,\;j=1,2\;,
\end{equation}
(we use the notation $\nu(.)$ defined in  \eqref{Gegen}). Thus, the general solution of \eqref{21} is
$$
\vec{y}(t,l)=\left(\,\,y_{1}(t,l),\,y_{2}(t,l),\,2ty'_{1}(t,l),\,2ty'_{2}(t,l)\,\,\right)^{Tr}.
$$

\subsection{Approximate general solution of the FD problem and further plan.}\label{subset:5}
Thus, in the regime
\begin{equation}\label{24}
\left\{\begin{array}{l}
         n\to\infty \\
         \displaystyle\frac{k}{n}\rightarrow t\in K\Subset(0,1]
       \end{array}
\right.\;,
\end{equation}
the general solutions (\ref{23}) could be a good approximation for
general solutions
$C_{1}\overrightarrow{Y}_{k}^{(1)}+C_{2}\overrightarrow{Y}_{k}^{(2)}$
of (\ref{17}) - (\ref{16}), for $k\in\mathbb{Z}$:
\begin{equation}\label{25}
\overrightarrow{Y}_{k}^{(1)}(\lambda)\approx\left(
\begin{array}{c}
y_{1}\left(\frac{k}{n},\lambda n^{4}\right)\\
0 \\
\frac{2k}{n} y'_{1}\left(\frac{k}{n},\lambda n^{4}\right) \\
0\\
\end{array}
\right)\;, \quad
\overrightarrow{Y}_{k}^{(2)}(\lambda)\approx\left(
\begin{array}{c}
0\\
y_{2}\left(\frac{k}{n},\lambda n^{4}\right)\\
0 \\
\frac{2k}{n} y'_{2}\left(\frac{k}{n},\lambda n^{4}\right)\\
\end{array}
\right)\;,
\end{equation}
here ($'$) denotes the derivative with respect to the first variable.


Now, we recall (see Introduction) the further steps we need to proceed in order to obtain in the regime
(\ref{24}) asymptotics of the exact constant (\ref{12}) and the
vector $\vec{v}$ - (\ref{14}) which defines the extremal polynomial
$Q'_{n}$. To choose from the general (approximate) solution
(\ref{25}) for $k\in\mathbb{Z}$, a solution which corresponds to the
boundary conditions (BC)
\begin{equation}\label{26}
v_{-2}=v_{-1}=0\quad\mbox{and}\quad v_{n}=v_{n+1}=0\;,
\end{equation}
we proceed like in \cite{Tul}. First we find a set of two particular solutions of (\ref{16})
$x_{k}(\lambda)$ for $\lambda=0$, which correspond to the boundary
conditions (\ref{26}) at the left end, i.e.
\begin{equation}\label{27}
x_{-2}^{(j)},\;x_{-1}^{(j)}=0\;,\quad j=1,2.
\end{equation}
The second step is to choose  constants $C_1,  C_2$ for (\ref{25}) such that
asymptotics of
$
y_{j}(t,l),\,\, j=1,2\;,$ when $t\to 0$
would match with asymptotics of
$
\overrightarrow{Y_{k}}^{(j)}(0)\,\,\,j=1,2\;$ when $k\to\infty$, which correspond to the particular
solutions $x_{k}^{(j)}(0)$.
It defines the unknown constants for (\ref{25}), (\ref{23}). The
last step is to obtain a linear combination of two approximate
discrete solutions $\widetilde{v}_{k}^{(j)}(\lambda),\;j=1,2$ (we
get them from (\ref{25})), satisfying the left end BC in
(\ref{26}), such that this combination satisfies the right end BC
in (\ref{26}). It is possible to do if
\begin{equation}\label{30}
\det\left|\begin{array}{ll}
\widetilde{v}_{n}^{(1)}(\lambda) & \widetilde{v}_{n+1}^{(1)}(\lambda) \\
\widetilde{v}_{n}^{(2)}(\lambda) & \widetilde{v}_{n+1}^{(2)}(\lambda) \\
\end{array}\right|=0\;.
\end{equation}
This equation (in $\lambda$) has the same meaning as (\ref{13}) and
its solutions with approximate $\widetilde{v}_{n}^{(j)}(\lambda)$
gives an approximation of the eigenvalues $\lambda(\mathbf{A},\mathbf{D})$.

\subsection{Two particular solutions of the FD problem for $\lambda=0$.}\label{subsec:6}
Here we find solutions of recurrence equation associated with \eqref{5-term} for $\lambda=0$, satisfying the BC (\ref{26}) at the left
end
\begin{equation}\label{32}
v_{-2}=v_{-1}=0\;.
\end{equation}
We are looking for the solutions of recurrences
\begin{equation}\label{31}
[\mathbf{C}_{1}^{Tr}\mathbf{C}_{2}^{Tr}\mathbf{N}^{-1}\mathbf{D}^{+}
\mathbf{N}^{-1}\mathbf{C}_{2}\mathbf{C}_{1}\vec{v}]_k\,\,= 0, \quad k=0,1,\ldots, n-3
\end{equation}
First we consider the equations
\begin{equation}\label{33}
[\mathbf{C}_{2}\vec{h}]_k=0\;,\ \mbox{with} \ (\mathbf{C}_{2}\vec{h})_{k-1}=h_{k-1}+
\frac{2k(k+\alpha+1)\,h_{k}}{(2k+\alpha+\beta+1)(2k+\alpha+\beta+2)}\;.
\end{equation}
The structure of this equation is such that $h_{-1}=0$ and $h_{-2}$
can take any values and we can put  it to be zero. It is easy to
check that the homogeneous equation (\ref{33}) has a solution
\begin{equation}\label{34}
h_{-2}=h_{-1}=0\,,\quad h_{k}=
\frac{(2k+\alpha+\beta+2)!}{(-2)^{k}k!(k+\alpha+1)!}\;,\quad k=0,1,\ldots,n-1\;.
\end{equation}
Next, we consider the equation:
\begin{equation}\label{35}
[\mathbf{C}_{1}\vec{h}]_k=0\;,\ \mbox{with} \ (\mathbf{C}_{1}\vec{h})_{k-1}=h_{k-1}+
\frac{-2k(k+\beta)}{(2k+\alpha+\beta)(2k+\alpha+\beta+1)}\,h_{k}\;.
\end{equation}
Again we can check that this homogeneous equation has a solution
$$
h_{-2}=h_{-1}=0\,,\quad h_{k}=
\frac{(2k+\alpha+\beta+1)!}{2^{k}k!(k+\beta)!}\;,\quad k=0,1,\ldots,n-1\;,
$$
which can be taken as the first particular solution of (\ref{31}),
(\ref{32}):
\begin{equation}\label{36}
v_{k}^{(1)}(0)=\frac{(2k+\alpha+\beta+1)!}{2^{k}k!(k+\beta)!}\;,\qquad k=0,1,\ldots,n-1\;.
\end{equation}
We find the second particular solution of (\ref{31}), (\ref{32}) as
a solution of the non homogeneous equation (\ref{35}) with right
hand side $h_{k-1}$ from (\ref{34})
$$
v_{k-1}^{(2)}(0)-\frac{2k(k+\beta)}{(2k+\alpha+\beta)(2k+\alpha+\beta+1)}\,v_{k}^{(2)}(0)=
\frac{-2k(2k+\alpha+\beta)!}{(-2)^{k}k!(k+\alpha)!}\;.
$$
It is easy to check that this non homogeneous equation has a
solution
\begin{equation}\label{37}
v_{k}^{(2)}(0)=\frac{(2k+\alpha+\beta+1)!}{(-2)^{k}k!(k+\alpha)!}\;,\qquad k=0,1,\ldots,n-1\;.
\end{equation}
Note, that for $v^{(1)}$, $v^{(2)}$ one has
\[
\mathbf{C}_{1}^{Tr}\mathbf{C}_{2}^{Tr}\mathbf{N}^{-1}\mathbf{D}^{+}
\mathbf{N}^{-1}\mathbf{C}_{2}\mathbf{C}_{1}\vec{v}\,\,=
\,\,\mbox{const}_1\, \vec{e}_{n-1} \,+ \,\mbox{const}_2 \vec{e}_{n-2}\;.
\]
where $e_{n-2}=(0,\ldots,0,1,0)^{Tr}$, $e_{n-1}=(0,\ldots,0,0,1)^{Tr}$.

\subsection{Matching  of the particular solutions of FD
and ODEs problems}\label{subset:7} We obtain the particular solutions
$\overrightarrow{X}_{k}^{(1)}(0),\overrightarrow{X}_{k}^{(2)}(0)$ of
(\ref{16})  from (\ref{36}), (\ref{37}),
 (\ref{15})
\begin{equation}\label{38}
x_{k}^{(1)}(0):=\frac{(k+\alpha)!}{k!}\;,\quad \quad x_{k}^{(2)}(0):=\frac{(k+\beta)!}{(-1)^{k}\,k!}\;.
\end{equation}
Using (\ref{17}) we have
\begin{equation}\label{39}
\overrightarrow{Y}_{k}^{(j)}(0)=
\left(
  \begin{array}{l}
    x_{k-2}^{(j)}+x_{k-1}^{(j)} \\
    x_{k-2}^{(j)}-x_{k-1}^{(j)} \\
    k\left(-(x_{k-2}^{(j)}+x_{k-1}^{(j)})+(x_{k}^{(j)}+x_{k+1}^{(j)})\right) \\
    k\left(-(x_{k-2}^{(j)}-x_{k-1}^{(j)})+(x_{k}^{(j)}-x_{k+1}^{(j)})\right) \\
  \end{array}
\right),\quad j=1,2\,.
\end{equation}
Substituting the expansions of (\ref{38}) in (\ref{39}),
we get for $\overrightarrow{Y}_{k}^{(1)}(0)$ and $\overrightarrow{Y}_{k}^{(2)}(0)$
$$
k^{\alpha}\left[
\left(
  \begin{array}{c}
    2+\frac{(\alpha-2)\alpha}{k} \\
    -\frac{\alpha}{k} \\
    4\alpha+O(\frac{1}{k}) \\
    O(\frac{1}{k}) \\
  \end{array}
\right)+O\left(\frac{1}{k^{2}}\right)\right]\;,
\quad
k^{\beta}\left[
\left(
  \begin{array}{c}
  -\frac{\beta}{k} \\
  2+\frac{(\beta-2)\beta}{k} \\
  O(\frac{1}{k}) \\
  4\beta+O(\frac{1}{k}) \\
  \end{array}
\right)+O\left(\frac{1}{k^{2}}\right)\right]\;,
$$
correspondingly. We conclude for $k\to\infty\,,\,\,\,\,j=1,2$ (recalling the notation \eqref{notbj}):
\begin{equation}\label{40}
\overrightarrow{Y}_{k}^{(j)}(0)=k^{b_{j}}\left(\overrightarrow{C}_{0}^{(j)}+O\left(\frac{1}{k}\right)\right),
\overrightarrow{C}_{0}^{(1)}:=\left(
                                  \begin{array}{c}
                                    2 \\
                                    0 \\
                                    4\alpha \\
                                    0 \\
                                  \end{array}
                                \right),
\overrightarrow{C}_{0}^{(2)}:=\left(
                                  \begin{array}{c}
                                    0 \\
                                    2 \\
                                    0 \\
                                    4\beta \\
                                  \end{array}
                                \right).
\end{equation}
Now, we state a "matching condition" for the choice of the
particular solutions of the differential problem (\ref{21})-(\ref{23}) when $t\to 0$:
\begin{equation}\label{41}
\vec{y}^{(j)}(t,l)=t^{b_{j}}\left(\overrightarrow{C}_{0}^{(j)}+o\left(1\right)\right)\;.
\end{equation}
If this condition is satisfied, then, for $\lambda=l\,n^{-4}$, we expect that in the regime
(\ref{24})
\begin{equation}\label{42}
\overrightarrow{Y}_{k}^{(j)}(l\,n^{-4})=n^{b_{j}}\vec{y}^{(j)}(\frac{k}{n},l)+o\left(k^{b_{j}}\right)\;.
\end{equation}
This assertion will be proved later (see Theorem~\ref{Th2}).\\
 Using the well known power series expansion  of the Bessel functions for $x\to 0$ (see \cite{BOW}),
$$
J_{\nu}(x)=\hat{c}_{\nu}x^{\nu}\left(1+\tilde{c}_{\nu}x^{2}+O(x^{4})\right)\,,\qquad \hat{c}_{\nu}:=
\frac{1}{2^\nu\;\nu!} ,
$$
and matching the conditions (\ref{41}), (\ref{40}), we obtain
expressions for the constants $\tilde{C}_{k,j}$ in the presentation
(\ref{23}) of the general solution of (\ref{21}):
$$
\tilde{C}_{2,j}=0\,,\quad \tilde{C}_{1,j}=2^{\nu_j+1}\left(\hat{c}_{\nu_{j}}l^{\nu_j/2}\right)^{-1}\,,\quad
\nu_{j}=\nu(b_{j})\,,\quad j=1,2.
$$
Thus, the particular solutions of (\ref{21}), satisfying  the
condition (\ref{41}), are
\begin{equation}\label{43}
\vec{y}^{(1)}(t,l)=\left(
                     \begin{array}{c}
                       y_{1}(t,l) \\
                       0 \\
                      2t \frac{d}{dt}\,y_{1}(t,l) \\
                       0 \\
                     \end{array}
                   \right),\qquad
\vec{y}^{(2)}(t,l)=\left(
                     \begin{array}{c}
                       0 \\
                       y_{2}(t,l) \\
                       0 \\
                       2t \frac{d}{dt}\,y_{2}(t,l) \\
                     \end{array}
                   \right),\end{equation}
where
\begin{equation}\label{44}
y_{j}(t,l)=\displaystyle\frac{2^{b_j}\nu_{j}!}{l^{\nu_{j}/2}}\,t\,J_{\nu_{j}}\left(\frac{\sqrt{l}\,t^{2}}{2}
\right)\,,\quad j=1,2.
\end{equation}

\subsection{Matching of the right end BC for $\{v_{k}^{(j)}\}_{k=0}^{n-1}$.}\label{8}
Now, substituting (\ref{43})-(\ref{44}) into (\ref{42}), (\ref{17}) and
(\ref{15}), we arrive to the two particular approximate sequences
$\{v_{k}^{(s)}(\lambda)\}_{k\in\mathbb{Z}_{+}}$, satisfying  the
left end BC in (\ref{26}):
$$
v_{n}^{(j)}(\lambda)\Bigr|_{\lambda=\frac{l}{n^{4}}}=\,(-1)^{(j-1)n}\,\frac{(2n+\alpha+\beta)!\,\,
\nu_{j}!\,\,n^{b_{j}}}
{2^{n-2\nu_{j}}\,(n+\alpha)!\, (n+\beta)! \, l^{\nu_{j}/2}}\,\left(J_{\nu_{j}}\left(\frac{\sqrt{l}}{2}\right)
+o\,(1)\right)\;.
$$
Thus, the right end BC (see (\ref{30})) when $n\to\infty$ and
$\lambda=\frac{l}{n^{4}}$, is equivalent to
$$
C(l,n)\left(\det\left|\begin{array}{cc}
J_{\nu(\alpha)}\left(\frac{\sqrt{l}}{2}\right) & J_{\nu(\alpha)}\left(\frac{\sqrt{l}}{2}\right) \\
-J_{\nu(\beta)}\left(\frac{\sqrt{l}}{2}\right) & J_{\nu(\beta)}\left(\frac{\sqrt{l}}{2}\right)
                 \end{array}
\right|\ +o(1)) \right)=0\;,\quad C(l,n)\neq 0\,,\;l>0\;.
$$
From here we conclude that the roots  of the equation
\begin{equation}\label{45}
J_{\nu(\alpha)}\left(\frac{\sqrt{l}}{2}\right)\,J_{\nu(\beta)}\left(\frac{\sqrt{l}}{2}\right)
+o(1)=0
\end{equation}
give approximate values of $\lambda=\frac{l}{n^{4}}$, for which
the BC (\ref{26}) are fulfilled, and the minimal root of (\ref{45}):
$$
\frac{\sqrt{l_{\alpha,\beta}}}{2}:= \min\limits_{\nu_{i}}\,(j_{\,\nu_{i}})=:j_{\,\nu^{*}}\;,
$$
(here $j_{\nu}$ is the minimal root of the Bessel function
$J_{\nu}$) (here we use a monotonicity on parameter $\nu$ of the minimal zero  of Bessel functions
$J_{\nu}$), see \cite{BOW}) gives the main term of asymptotics (see (\ref{12}) and
(\ref{13})) of the exact Markov-Bernstein constant
\begin{equation}\label{46}
M_{n}=\frac{n^{2}}{2j_{\,\nu^{*}}}\,(1+o(1))\;.
\end{equation}

\section{Matching and convergence of  FD and  DE problems}\label{sec:3}
\subsection{Statements of the results.}\label{9}
Our derivation of the asymptotics (\ref{46}) contained one assumption
which requires a special rigorous treatment. It is the convergence
in the regime \eqref{24} of the discrete solution to the continuous solution,
see (\ref{42}). Here we state a theorem\ which establishes (\ref{42})
under a restriction on $(\alpha,\beta)$.

\begin{theorem}\label{Th2} Let
$\{\overrightarrow{Y}_{k}^{(j)}(\lambda)\}_{k=0}^{\infty},\,j=1,2$,
be a set of the particular solutions of the FD problem, i.e. the recurrence \eqref{17a}
with  the   matrix
$\widehat{\mathbf{M}}_{2}^{(\alpha,\beta)}(k,\lambda)$ such that
$\{\overrightarrow{Y}_{k}^{(j)}(\lambda)\}_{k=0}^{\infty},\;j=1,2$,
are given
 by
\eqref{39}, \eqref{38}, $\overrightarrow{Y}_{0}^{(j)}(\lambda)=\overrightarrow{Y}_{0}^{(j)}(0)$

Let the parameters $(\alpha,\beta)$ in
$\widehat{\mathbf{M}}_{2}^{(\alpha,\beta)}$ satisfy the condition:
\begin{equation}\label{47}
 |\alpha - \beta| < 4 \,, \qquad\alpha,\beta > -1\;.
\end{equation}
Then, for $\lambda=\frac{l}{n^{4}},\,\frac{k}{n}\to t$ and
$n\to\infty$, uniformly for $l\in
\widetilde{K}\Subset\mathbb{C},\,t\in K \Subset (0,1]$,
\begin{equation}\label{48}
\overrightarrow{Y}_{k}^{(j)}(\lambda)=n^{b_{j}}\left( \vec{y}^{(j)}(t,l)+o(t^{b_j})\right)\;,
\quad b_{j}=\left\{\begin{array}{cc}
                     \alpha, & j=1 \\
                     \beta, & j=2
                   \end{array}\right.\;,
\end{equation}
holds true. Here
$\vec{y}^{(j)},\,j=1,2,$ are the particular solutions  \eqref{44}, \eqref{43} of the
DE problem \eqref{21} satisfying the matching condition
\eqref{41}.
\end{theorem}
Taking
into account the procedure of derivation of (\ref{46}) we obtain the validity of  Theorem~\ref{T1} as a  corollary of  Theorem~\ref{Th2}.

\subsection{Proof of Theorem~\ref{Th2}.}\label{subsec:10}

 In our proof we use an approach proposed in ~\cite{Tul}.
The comparison of $\overrightarrow{Y}_{k}^{j}$ and
$n^{b_{j}}\vec{y}^{(j)}$, in the regime (\ref{24}), will be performed
by means of the following relations. In what follows we denote by $|v|$ a norm of the vector $v$ in the vector space $R^n$ and by $||S||$ the associated matrix norm. For two recurrent sequences
$\{\vec{v}_{k}\}$, $\{\vec{w}_{k}\}$, defined  by
$$
\vec{v}_{k+2}=S_{k}^{(1)}\vec{v}_{k}\,,\quad\vec{w}_{k+2}=S_{k}^{(2)}\vec{w}_{k}\,,\quad
\vec{v}_{k}, \vec{w}_{k}\in\mathbb{R}^{N},\quad k\in\mathbb{N}\ast 2 \;,
$$
we have for $k>k_{0}$
\begin{equation}\label{49}
 \begin{array}{l} \qquad \qquad \qquad
|\vec{w}_{k}-\vec{v}_{k}|\leqslant\\\left[|\vec{w}_{k_{0}}-\vec{v}_{k_{0}}|\!\!+\!\!
\sum\limits_{m=k_{0},2}^{k-2}\left|(S_{m}^{(2)}-S_{m}^{(1)})\vec{v}_{m}\right|E(k_{0},m,S^{(2)})\!\right]
E^{-1}(k_{0},k-2,S^{(2)}),
\end{array}
\end{equation}
where the summation is performed with the step 2, and
$$
E(k_{0},m,S):=\exp\left(\sum\limits_{i=k_{0},\;2}^{m}(1-\|S_i\|)\right)\;.
$$
The estimate (\ref{49}) easily follows by induction from
$$
|\vec{w}_{n+2}-\vec{v}_{n+2}|=\left|S_{n}^{(2)}\vec{w}_{n}-S_{n}^{(1)}\vec{v}_{n}\right|\leqslant
\|S_{n}^{(2)}\|\,|\vec{w}_{n}-\vec{v}_{n}|+\left|(S_{n}^{(2)}-S_{n}^{(1)})\vec{v}_{n}\right|
$$
and from the inequality $xe^{1-x}\leqslant 1$ applied to
$\|S_{n}^{(2)}\|$.

1) We are going to use the estimate (\ref{49}) to compare the
solution $\overrightarrow{Y}_{k}^{(j)}(\lambda)$ of the finite
difference problem (\ref{18})  with the solution
$\vec{y}^{(j)}(t,l)$ (see (\ref{43})-(\ref{44})) of the
differential problem (\ref{21}). To do this (having in mind
(\ref{49})) we define a difference operator which connects values
of $\vec{y}^{(j)}(t)$ taken on the discrete grid
$t:=\frac{k}{n},\,k\in\mathbb{N}_{+}$:
\begin{equation}\label{50}
\vec{y}^{(j)}\left(\frac{k+2}{n},\,l\right)=\widetilde{\widetilde{\mathbf{M}}}_{2}(k,l,n)\,\,
\vec{y}^{(j)}\left(\frac{k}{n},\,l\right),\quad \widetilde{\widetilde{\mathbf{M}}}_{2}=:
\left(\mathbf{I}+\frac{2}{k}\widetilde{\widetilde{\mathbf{M}}}_{3}\right)\,.
\end{equation}
We apply the following lemma from ~\cite{Tul} (see Lemma 3.2).

\begin{lemma}\label{L2} Let $M(t)$ be a matrix-valued function solving the
Cauchy problem with smooth matrix-valued coefficient $F$:
$$
\frac{d}{dt}M(t)=F(t)\,M(t),\qquad M(t_{1})=\mathbf{I}\,.
$$
Then the following estimate holds for
$t_{2}>t_{1}\,,\;\displaystyle\int\limits_{t_{1}}^{t_{2}}\|F(t)\|\,dt<1$:
$$
\|M(t_{2})-\mathbf{I}-(t_{2}-t_{1})\,F(t_{1})\|<2(t_{2}-t_{1})^{2}\max\limits_{[t_{1};\,t_{2}]}\,\|F'(t)+F^{2}(t)\|\,.
$$
\end{lemma}
Applying this lemma to $M=\widetilde{\widetilde{\mathbf{M}}}_{2}$. $t_1=\frac{k}{n}$, $t_2=\frac{k+2}{n}$  we get an estimate which we shall use in our analysis later.
\begin{equation}\label{51}
\left\|\widetilde{\widetilde{\mathbf{M}}}_{3}(k,l,n)-\widetilde{\mathbf{M}}_{3}\left(\frac{k}{n},l\right)\right\|=
O\left(\frac{1}{k}\right)\;,
\end{equation}
here $\widetilde{\mathbf{M}}_{3}$ is defined in (\ref{21}) and
$k\to\infty,\;\frac{k}{n}\in K_{1}\Subset\mathbb{R}_{+}\,,\;l\in
K_{2}\Subset\mathbb{C}$\;.

2) Now, we can rewrite  (\ref{49}) for our purpose
\begin{equation}\label{52}
\begin{array}{l}
\left|\overrightarrow{Y}_{k}^{(j)}(\lambda)-n^{b}\vec{y}^{(j)}
\left(\frac{k}{n},\lambda n^{4}\right)\right|\lesssim
\left\{\left|\overrightarrow{Y}_{k_{0}}^{(j)}(\lambda)-n^{b}\vec{y}^{(j)}
\left(\frac{k_{0}}{n},\lambda n^{4}\right)\right|\right.*\\
\displaystyle\sum\limits_{m=k_{0},2}^{k-2}\left|\left(\widehat{\mathbf{M}}_{2}^{(\alpha,\beta)}(m,\lambda)-
\widetilde{\widetilde{\mathbf{M}}}_{2}(m,\lambda n^{4})\right)\,n^{b}\vec{y}^{(j)}
\left(\frac{k}{n},\lambda n^{4}\right)\right|*\\
\left.E\left(k_{0},m,\widehat{\mathbf{M}}_{2}^{(\alpha,\beta)}\right)
\right\}E^{-1}\left(k_{0},k-2,\widehat{\mathbf{M}}_{2}^{(\alpha,\beta)}\right)\;,\\
\end{array}
\end{equation}
where
$\widehat{\mathbf{M}}_{2},\widetilde{\widetilde{\mathbf{M}}}_{2}$ are defined in
(\ref{17a}), (\ref{50}) correspondingly.

3) To proceed with (\ref{52}) we start with an estimation of the
initial deviation. We put $k_{0}:=n^{1-\delta}$, where $\delta>0$
will be fixed later. We shall compare the initial data in
(\ref{52}), estimating their deviation from the corresponding
particular solution $Y_{k_{0}}^{(j)}(0)$, which we know explicitly
 (\ref{39}), (\ref{38})
\begin{equation}\label{53}
\begin{array}{l}
\left|\overrightarrow{Y}_{k_{0}}^{(j)}(\lambda)-n^{b}\vec{y}^{(j)}
\left(\frac{k_{0}}{n},\lambda n^{4}\right)\right|< \\
\\
\left|\overrightarrow{Y}_{k_{0}}^{(j)}(\lambda)-
\overrightarrow{Y}_{k_{0}}^{(j)}(0)\right|+
\left|\overrightarrow{Y}_{k_{0}}^{(j)}(0)-n^{b}\vec{y}^{(j)}
\left(\frac{k_{0}}{n},\lambda
n^{4}\right)\right|\,,\quad k_{0}=n^{1-\delta}\,.
\end{array}
\end{equation}
To estimate the second term in the right hand side of (\ref{53}) we
use the matching condition (\ref{41})
\begin{equation}\label{54}
\begin{array}{l}
\left|\overrightarrow{Y}_{k_{0}}^{(j)}(0)-n^{b_j}\vec{y}^{(j)}
\left(\frac{k_{0}}{n},l\right)\right|=\qquad \\
\\
\qquad  \qquad \qquad \qquad \left|\left(\overrightarrow{C}_{0}^{(j)}k_{0}^{b_j}
+O\left(\frac{1}{k_{0}}\right)\right)-
\left(\overrightarrow{C}_{0}^{(j)}k_{0}^{b_j}+O\left(\frac{k_{0}^{4}}{n^{4}}\right)\right)\right|\,=\\\\
\qquad \qquad \qquad \qquad  \qquad \qquad \qquad \qquad \qquad
O\left(k_{0}^{b-1}\right)+O\left(k_{0}^{b}\left(\frac{k_{0}}{n}\right)^{4}\right).
\end{array}\end{equation}
To estimate the first term in the right hand side of (\ref{53}) we
again use (\ref{49})
\begin{equation}\label{55}
\begin{array}{l}
\left|\overrightarrow{Y}_{k_{0}}^{(j)}(\lambda)-
\overrightarrow{Y}_{k_{0}}^{(j)}(0)\right|\leqslant \left\{\left|\overrightarrow{Y}_{2}^{(j)}(\lambda)-
\overrightarrow{Y}_{2}^{(j)}(0)\right|*\right.\\
\qquad \qquad \displaystyle\sum\limits_{m=2,2}^{k_{0}-2}
\left|\left(\widehat{\mathbf{M}}_{2}^{(\alpha,\beta)}(m,\lambda)-
\widehat{\mathbf{M}}_{2}^{(\alpha,\beta)}(m,0)\right)\overrightarrow{Y}_{m}^{(j)}(0)\right|*\\
\qquad \qquad \qquad \qquad \left.E\left(k_{0},m,\widehat{\mathbf{M}}_{2}^{(\alpha,\beta)}(m,\lambda)\right)\right\}
E^{-1}\left(k_{0},m,\widehat{\mathbf{M}}_{2}^{(\alpha,\beta)}(m,\lambda)\right).
\end{array}
\end{equation}
For further estimations we shall use $\|\quad\|_{\delta}$-norm
introduced in ~\cite{Tul}. This norm is related with a basis in
which the operator $\mathbf{M}_{3}^{(0)}(\infty)$ defined in  (\ref{17a}), (\ref{18})
has a matrix (with eigenvalues on the diagonal) which is arbitrary
close to a diagonal matrix
$$
\lim\limits_{m\to\infty}\widehat{\mathbf{M}}_{2}^{(\alpha,\beta)}(m,0)=:\widehat{\mathbf{M}}_{2}^{(\alpha,\beta)}(\infty,0)=
\left(\mathbf{I}+\frac{1}{k}\mathbf{M}_{3}^{(0)}(\infty)\right)\,.
$$
In such a basis the operator
$\widehat{\mathbf{M}}_{2}^{(\alpha,\beta)}(\infty,0)$ will have a norm close to $(1+\frac{1}{k}\,B)$,
where
\begin{equation}\label{57}
 B:=\max\{\alpha,2-\alpha,\beta,2-\beta\}\
 \end{equation}
 is the maximal real part of eigenvalues of the matrix
$$ \mathbf{M}_{3}^{(0)}(\infty):=\widetilde{\mathbf{M}}_{3}(0,0)=
 \left(
  \begin{array}{cccc}
    0 & 0 & 1/2 & 0 \\
    0 & 0 & 0 & 1/2 \\
    2\alpha(\alpha-2) & 0 & 2 & 0 \\
    0 & 2\beta(\beta-2) & 0 & 2 \\
    \\
  \end{array}\right).
$$
That is, for each $\delta>0$ there exist $k_{\delta}\in\mathbb{N}$
and a basis such that, with respect to the Euclidean norm
$|\quad|_{\delta}$ associated with this basis (see \cite{Tul} for the details about $\delta$-norm), the corresponding
operator norm $\|\quad\|_{\delta}$ has the estimate
\begin{equation}\label{58}
\forall\, k>k_{\delta},\qquad\left\|\widehat{\mathbf{M}}_{2}^{(\alpha,\beta)}(k,0)\right\|=
\left\|\mathbf{I}+\frac{1}{k}\mathbf{M}_{3}^{(0)}(k)\right\|_{\delta}<1+\frac{B+\delta}{k}\;.
\end{equation}
Then, we estimate the terms from the right hand side of (\ref{55}).
We have (see (\ref{18}))
$$\begin{array}{c}
\left|\left(\widehat{\mathbf{M}}_{2}^{(\alpha,\beta)}(m,\lambda)-
\widehat{\mathbf{M}}_{2}^{(\alpha,\beta)}(m,0)\right)\overrightarrow{Y}_{m}^{(j)}(0)\right|_{\delta}=\\
\\
=\frac{1}{m}\left|\left(\mathbf{M}_{3}(m,\lambda)-
\mathbf{M}_{3}(m,0)\right)\overrightarrow{Y}_{m}^{(j)}(0)\right|_{\delta}=
\frac{1}{m}\left|\mathbf{M}_{3}^{(1)}\lambda\overrightarrow{Y}_{m}^{(j)}(0)\right|\leqslant
O\left(|\lambda|\frac{m^{4}}{m}m^{b_{j}}\right)\,.\end{array}
$$
Analogously
$$
\|\mathbb{\mathbf{M}}_{3}(m,\lambda)-\mathbf{M}_{3}(m,0)\|_{\delta}=|\lambda| \|\mathbf{M}_{3}^{(1)}\|=O(\lambda m^{4})\;,
$$
and
$$
\left\|\;\mathbf{I}+\frac{1}{m}\mathbf{M}_{3}(m,\lambda)\right\|_{\delta}\leqslant1+\frac{B+\delta}{m}+
\frac{1}{m}\left\|\mathbf{M}_{3}(m,\lambda)-\mathbf{M}_{3}(m,0)\right\|_{\delta}=
1+\frac{B+\delta+O(|\lambda|m^{4})}{m}\,.
$$
Now, we proceed with (\ref{55}). For $k=k_{0}=n^{1-\delta}$ we have
$$
\begin{array}{c}
\left|\overrightarrow{Y}_{k}^{(j)}(\lambda)-\overrightarrow{Y}_{k}^{(j)}(0)\right|\leqslant O(1)
\exp\left\{\sum\limits_{\tilde{k}=2,2}^{k-2}\frac{B+\delta+O(|l|n^{-4})\tilde{k}^{4}}{\tilde{k}}\right\}*
\qquad \qquad \qquad \qquad \qquad\\
\\
\qquad \qquad \qquad \sum\limits_{\tilde{k}=2,2}^{k-2}\frac{O(|l|n^{-4})\tilde{k}^{4+b_{j}}}{k}
\exp\left\{-\sum\limits_{\tilde{\tilde{k}}=2,2}^{\tilde{k}-2}\frac{B+\delta+O(|l|n^{-4})\tilde{\tilde{k}}^{4}}
{\tilde{\tilde{k}}}\right\}=\\
\\
\frac{|l|}{n^{4}}\,O\left(k^{B+\delta}\exp\left[O\left(|l|\left(\frac{k}{n}\right)^{4}
\right)\right]\right)\,\sum\limits_{\tilde{k}=2,2}^{k-2}\tilde{k}^{4+b_{j}-B-\delta-1}
\exp\left[O\left(|l|\left(\frac{\tilde{k}}{n}\right)^{4}\right)\right].
\end{array}
$$
At this point we assume, that
$$
4+b_{j}>B+\delta\;.
$$
This assumption (see \eqref{57}) implies the restriction (\ref{47}) in our theorems. Thus,
we obtained
$$
\left|\overrightarrow{Y}_{k_{0}}^{(j)}(\lambda)-\overrightarrow{Y}_{k_{0}}^{(j)}(0)\right|\leqslant O\left(
\left(\frac{k_{0}}{n}\right)^{4}k_{0}^{b_{j}}\right)\;,
$$
and finally, for the initial deviation (\ref{53}) in (\ref{52}) we
get from here and (\ref{54})
\begin{equation}\label{59}
\left|\overrightarrow{Y}_{k_{0}}^{(j)}(\lambda)-n^{b_{j}}\overrightarrow{y}
\left(\frac{k_{0}}{n},\lambda n^{4}\right)\right|\leqslant O\left(k_{0}^{b_{j}-1}\right)+
O\left(\left(\frac{k_{0}}{n}\right)^{4}k_{0}^{b_{j}}\right)\;.
\end{equation}

4) Now, we come back to (\ref{52}). Using (\ref{51}) and triangle inequality one has
$$
\begin{array}{c}
k\left\|\;\widetilde{\widetilde{\mathbf{M}}}_{2}(k,\lambda n^{4})-\widehat{\mathbf{M}}_{2}^{(\alpha,\beta)}
(k,\lambda)\right\|_{\delta}=\left\|\;\widetilde{\widetilde{\mathbf{M}}}_{3}(k,l,n)-\mathbf{M}_{3}\left(k,\frac{l}{n^{4}}\right)
\right\|_{\delta}\leqslant\\
\\
O\left(\frac{1}{k}\right)+\frac{l}{n^{4}}(k^{4}+\mbox{const}\,k^{3}+\ldots)-\frac{l}{n^{4}}k^{4}=
O\left(\frac{1}{k}\right)+O\left(\frac{l}{n^{4}}k^{3}\right)=\\
\\
o\left(\left(\frac{k}{n}\right)^{4}k^{-\varepsilon_{1}}\right)+o\left(k^{-\varepsilon}\right)\,,
\qquad \qquad \qquad
\forall\,\varepsilon,\varepsilon_{1}\in(0,1)\;. \qquad
\end{array}
$$
Substituting this estimate and (\ref{59}), (\ref{58}) in (\ref{52}),
we proceed
$$
\begin{array}{c}
  \left|n^{b_{j}}\vec{y}^{(j)}\left(\frac{k}{n},\lambda n^{4}\right)-
  \overrightarrow{Y}_{k}^{(j)}(\lambda)\right|\leqslant\Biggl\{O\left(k_{0}^{b_{j}-1}\right)+
  O\left(k_{0}^{b_{j}}\left(\frac{k_{0}}{n}\right)^{4}\right)+  \\
  \\
  \sum\limits_{m=k_{0},2}^{k-2}\frac{1}{m}\left[o(m^{-\varepsilon_{1}}\left(\frac{m}{n}\right)^{4})
+o(m^{-\varepsilon})\right]n^{b_{j}}O\left(\frac{m}{n}\right)^{b_{j}}
\exp\left(-\sum\limits_{s=k_{0},2}^{m}\frac{B+\delta+O(|l|n^{-4}s^{4})}{s}\right)\Biggr\}* \\
  \\
  \exp\left(\sum\limits_{s=k_{0},2}^{m}\frac{B+\delta+O(|l|n^{-4}s^{4})}{s}\right)\;.
\end{array}
$$
Assuming again from (\ref{47})
$$
4+b_{j}>B+\delta+\varepsilon_{1}\;,
$$
we continue
$$
\begin{array}{c}
  \left|n^{b}\vec{y}^{(j)}\left(\frac{k}{n},\lambda n^{4}\right)-
  \overrightarrow{Y}_{k}^{(j)}(\lambda)\right|\leqslant\Biggl\{o\left(k_{0}^{b_{j}-\tilde{\varepsilon}}
  \right)+o\left(k_{0}^{b_{j}}\left(\frac{k_{0}}{n}\right)^{B+\delta-b_{j}}\right)+  \\
  \\
  \sum\limits_{m=k_{0},2}^{k}\frac{1}{m}\left[o(m^{b-\varepsilon_{1}}\left(\frac{m}{n}\right)^{4})
+o(m^{b-\varepsilon})\right]*\left(\frac{k_{0}}{m}\right)^{B+\delta}
\exp\left(O\left[\left(\frac{m}{n}\right)^{4}-\left(\frac{k_{0}}{n}\right)^{4}\right]\right)\Biggr\}* \\
  \\
  \left(\frac{k}{k_{0}}\right)^{B+\delta}
\exp\left(O\left[\left(\frac{k}{n}\right)^{4}-\left(\frac{k_{0}}{n}\right)^{4}\right]\right)\leqslant
\Bigl[o\left(k_{0}^{b_{j}-\tilde{\varepsilon}}\right)\left(\frac{k}{k_{0}}\right)^{B+\delta}+o(n^{b_{j}})+\\
\\
  \sum\limits_{m=k_{0},2}^{k-2}\frac{1}{m}\left[o\left(m^{b_{j}-\varepsilon_{1}}\left(
  \frac{m}{n}\right)^{4}\right)+o\left(m^{b_{j}-\varepsilon}\right)\right]
  \left(\frac{k}{m}\right)^{B+\delta}\exp(O(1))\Bigr]\exp(O(1))\leqslant \\
  \\
 o(n^{b})+o\left(k_{0}^{b_{j}-\tilde{\varepsilon}}\left(\frac{k}{k_{0}}\right)^{B+\delta}\right)+
  o\left(k^{b_{j}-\varepsilon_{1}}\right)+o\left(k_{0}^{b_{j}-\tilde{\varepsilon}}
  \left(\frac{k}{k_{0}}\right)^{B+\delta}\right)\;.
\end{array}
$$
Taking into account that
$$
k_{0}=n^{1-\delta}\,,\qquad\frac{k}{n}=t\in K\Subset\mathbb{R}_{+}\;,
$$
we see that the right hand side of this estimation is equal to
$$
o\left(n^{\max(\,b_{j},\,\, b_{j}-\tilde{\varepsilon}+(B-b_{j}+\tilde{\varepsilon}+\delta)\delta,\;\,
b_{j}-\varepsilon_{1},\;\, b_{j}-\varepsilon+(B-b_{j}+\varepsilon+\delta)\delta\,)}\right)\;.
$$
Since $\tilde{\varepsilon}<\varepsilon$, we choose $\delta$:
$$
\delta(B-b+\tilde{\varepsilon}+\delta)<\tilde{\varepsilon}\;,
$$
which yields
$$
\left|n^{b_{j}}\vec{y}^{(j)}\left(\frac{k}{n},l\right)-\overrightarrow{Y}_{k}^{(j)}
\left(\frac{l}{n^{4}}\right)\right|=o(n^{b_{j}})\;.
$$
Theorem is proved.

\end{document}